\documentclass[12pt]{amsart}
\advance\hoffset by -0.5 truecm
\usepackage{amssymb}
\newtheorem{Theorem}{Theorem}[section]
\newtheorem{Lemma}[Theorem]{Lemma}

\newtheorem{Remark}[Theorem]{Remark}

\def \dim{{\mbox {dim}}\,}

\def\V{\mbox{Var}}

\def\Div{\mbox{div}}
\def\Z{{\mathbb Z}}
\def\R\re
\def\V{\bf V}

\def \la{\lambda}

\def \re{{\mathbb R}}

\def \K{{\mathbb K}}
\def \V{{\bf V}}

\def \0{\lambda_{0}}
\def \la{\lambda}
\def \ga{\gamma}

\begin{document}
\title[Entropy production in Gaussian thermostats]{Entropy production in Gaussian thermostats}

\author[N.S. Dairbekov]{Nurlan S. Dairbekov}
\address{Kazakh British Technical University,
Tole bi 59, 050000 Almaty, Kazakhstan }
\email{Nurlan.Dairbekov@gmail.com}

\author[G.P. Paternain]{Gabriel P. Paternain}
 \address{ Department of Pure Mathematics and Mathematical Statistics,
University of Cambridge,
Cambridge CB3 0WB, England}
 \email {g.p.paternain@dpmms.cam.ac.uk}

%\subjclass{53C25, 53C21, 58F17, 35J15}

%\date{January 2006}

%\maketitle

\begin{abstract} We show that an arbitrary Anosov Gaussian thermostat on a surface is dissipative
unless the external field has a global potential.
This result is obtained by studying the cohomological equation of more general thermostats
using the methods in \cite{DP1}.

\end{abstract}

\maketitle

\section{Introduction}

Gaussian thermostats provide interesting models in nonequilibrium statistical mechanics
\cite{Ga,GaRu,Ru1}. Given a closed Riemannian manifold $(M,g)$ and a vector field $E$ (the {\it external
field}) on $M$,
the {\it Gaussian thermostat} (or isokinetic dynamics, cf. \cite{H}) is given by the differential equation
\begin{equation}
\frac{D\dot{\gamma}}{dt}=E(\ga)-\frac{\langle E(\ga),\dot{\ga}\rangle}{|\dot{\ga}|^2}\,\dot{\ga}.
\label{eqt}
\end{equation}
This equation defines a flow $\phi$ on the unit sphere bundle $SM$ of $M$ which reduces to the geodesic
flow when $E=0$.

In general, Gaussian thermostats are not volume preserving and the purpose of the present paper is to characterize
precisely those Anosov Gaussian thermostats in $2$ degrees of freedom which do not preserve
any smooth measure.

When $\phi$ is Anosov and $\dim M=2$ a result of E. Ghys \cite{Ghy} ensures that $\phi$ is topologically
conjugate to the geodesic flow of a metric of constant negative curvature and thus $\phi$ is transitive and 
topologically mixing.
For such a flow it is well known (cf. \cite[Chapter 20]{KH}) that there exists a unique Gibbs state $\rho$ associated with the
H\"older continuous potential $-\left.\frac{d}{dt}\right|_{t=0}\log J_t^u$, where $J_t^u$ is the {\it unstable Jacobian} of $\phi$.
The measure $\rho$ is characterized by being the maximum of
\[\nu\mapsto h_{\nu}(\phi)-\int\left.\frac{d}{dt}\right|_{t=0}\log J_t^u\,d\nu\]
where $\nu$ runs over all $\phi$-invariant Borel probability measures and $h_{\nu}(\phi)$ is the
measure theoretic entropy of $\phi$ with respect to $\nu$. The unique measure $\rho$ is called
the {\it SRB measure} of $\phi$. If $\tau$ is a probability measure which is absolutely continuous
with respect to the Liouville measure of $SM$, then $\rho$ is also the weak limit
of $\frac{1}{T}\int_{0}^{T}\phi^*\tau\,dt$ as $T\to\infty$.

{\it The entropy production} of the state $\rho$ is given by (cf. \cite{Ru0})
\[e_{\phi}(\rho):=-\int \Div F\,d\rho=-\sum\,\mbox{\rm Lyapunov exponents}\]
where $F$ is the infinitesimal generator of $\phi$ and $\Div F$ is the divergence of $F$ with respect
to any volume form in $SM$. 

Fix a volume form $\Theta$ on $SM$. Any other volume form can be written as $f\Theta$ for some 
smooth positive function $f$. 
If we let $L_F\Theta$ be the Lie derivative of $\Theta$ along $F$, then
\[L_{F}(f\Theta)=d(i_{F}f\Theta)=F(f)\Theta+fL_{F}\Theta=F(f)\Theta+f\Div F\Theta.\]
Hence if $\widetilde{\Div}F$ denotes the divergence of $F$ with respect to $f\Theta$ we have
\begin{equation}
\widetilde{\Div} F=F(\log f)+\Div F.
\label{coho}
\end{equation}
In other words the two divergences are flow cohomologous (and thus $e_{\phi}$ is well defined for any
$\phi$-invariant measure).

Ruelle \cite{Ru0} has shown that $e_{\phi}(\rho)\geq 0$
with equality if and only if $\rho$ is also the SRB measure of the flow $\phi_{-t}$.
If $\rho$ is an SRB measure for both $\phi_t$ and $\phi_{-t}$
then the theory of Gibbs states for Anosov flows (cf. \cite[Proposition 20.3.10]{KH})
implies that $-\left.\frac{d}{dt}\right|_{t=0}\log J_t^u$ and $\left.\frac{d}{dt}\right|_{t=0}\log J_t^s$
are cohomologous (and the coboundary is the derivative along the flow of a H\"older continuous function). 
It follows that $\phi$ preserves an absolutely continuous invariant measure with positive continuous density
(and this measure would have to be $\rho$). 
An application of the smooth Liv\v sic theorem 
\cite[Corollary 2.1]{LMM}
shows that $\phi$ preserves an absolutely continuous invariant measure with positive continuous density
if and only if $\phi$ preserves a smooth volume form. 
Using (\ref{coho}) we see that $e_{\phi}(\rho)=0$ if and only if $\Div F$ is a flow coboundary and we can
take $\Div F$ with respect to any volume form.

Let $\theta$ be the 1-form dual to $E$, i.e., $\theta_{x}(v)=\langle E(x),v\rangle$. An easy calculation
(see Lemma \ref{lied}) shows that if we consider in $SM$ the volume form determined by the canonical contact 1-form, then
$\Div F(x,v)=-\theta_{x}(v)$. Thus $e_{\phi}(\rho)=0$ if and only if there is a smooth solution
$u$ to the cohomological equation
\begin{equation}
F(u)=\theta.
\label{amiga}
\end{equation}
We will show as a consequence of a more general result to be stated below that if $\dim M=2$
then (\ref{amiga}) holds if and only if $\theta$ is an exact form, i.e. if and only if 
$E$ has a global potential. Thus we obtain:

\medskip

\noindent {\bf Theorem A.} {\it An Anosov Gaussian thermostat on a closed surface has
zero entropy production if and only if the external field $E$ has a global potential.
}

\medskip

A system with $e_{\phi}(\rho)>0$ is referred to as {\it dissipative}. Dissipative Gaussian thermostats
provide a large class of examples to which one can apply the Fluctuation Theorem of G. Gallavotti
and E.G.D. Cohen \cite{GC,GC1,Ga0} (extended to Anosov flows by G. Gentile \cite{G}) and this theorem is perhaps
one of the main motivations for determining precisely which thermostats are dissipative.
Observe that Gaussian thermostats are reversible in the sense that the flip
$(x,v)\mapsto (x,-v)$ conjugates $\phi_{t}$ with $\phi_{-t}$ (just as in the case of geodesic flows).
We recall that the {\it chaotic hypothesis} of Gallavotti and Cohen asserts that for systems
out of equilibrium, physically correct macroscopic results
will be obtained by assuming that the microscopic dynamics is uniformly hyperbolic.

In \cite{W1}, M. Wojtkowski proved Theorem A assuming that $E$ has a local potential (i.e. $\theta$
is closed) and in \cite{BGM}, F. Bonetto, G. Gentile and V. Mastropietro proved the theorem for the case
of a metric of constant negative curvature and $\theta$ a harmonic 1-form. 
We emphasize that we do not make any assumptions on $g$ or $E$ except that the underlying isokinetic dynamics
is Anosov. Conditions under which the Anosov property holds have been given in \cite{W1,W2}.

We now explain for which Anosov systems we can understand the cohomological equation (\ref{amiga}) completely.

Let $M$ be a closed manifold endowed with a Riemannian metric $g$.
We consider a {\it generalized isokinetic thermostat}. This consists of a semibasic vector field
$E(x,v)$, that is, a smooth map $TM \ni (x,v)\mapsto E(x,v)\in TM$ such that $E(x,v)\in T_{x}M$
for all $(x,v)\in TM$. As before the equation
\[\frac{D\dot{\gamma}}{dt}=E(\ga,\dot{\ga})-\frac{\langle E(\ga,\dot{\ga}),\dot{\ga}\rangle}{|\dot{\ga}|^2}\,\dot{\ga}.\]
defines a flow $\phi$ on the unit sphere bundle $SM$. These generalized thermostats
are no longer reversible unless $E(x,v)=E(x,-v)$.

Suppose now that $M$ is a closed oriented surface. We can write
\[E(x,v)=\kappa(x,v)v+\la(x,v)iv\]
where $i$ indicates rotation by $\pi/2$ according to the orientation
of the surface and $\kappa$ and $\la$ are smooth functions. The evolution of the thermostat on $SM$
can now be written as
\begin{equation}
\frac{D\dot{\gamma}}{dt}=\la(\gamma,\dot{\ga})\,i\dot{\ga}.
\label{eqgt}
\end{equation}
If $\la$ does not depend on $v$, then $\phi$ is the
{\it magnetic flow} associated with the magnetic field $\la\Omega_{a}$, where $\Omega_{a}$ is the area
form of $M$. Of course, magnetic flows are Hamiltonian.
If $\la$ depends linearly on $v$, we obtain the Gaussian thermostat (\ref{eqt}).

%In the present paper we shall study rigidity properties of {\it Anosov}
%magnetic flows. The Anosov property means that $T(SM)$
%splits as $T(SM)=E^{0}\oplus E^{u}\oplus E^{s}$ in such a way that
%there are constants $C>0$ and $0<\rho<1<\eta$ such that $E^{0}$ is
%spanned by  the generating vector field of the flow, 
%and for all $t>0$ we have
%\[\|d\phi_{-t}|_{E^{u}}\|\leq C\,\eta^{-t}\;\;\;\;\mbox{\rm
%and}\;\;\;\|d\phi_{t}|_{E^{s}}\|\leq C\,\rho^{t}.\]
%The subbundles are then invariant and H\"older continuous and have
%smooth integral manifolds, the stable and unstable manifolds,
%which define a continuous foliation with smooth leaves.

Let $\pi:SM\to M$ be the canonical projection.

\medskip

\noindent {\bf Theorem B.} {\it Let $M$ be a closed oriented surface and consider
a generalized isokinetic thermostat (\ref{eqgt}). Suppose the flow $\phi$ 
is Anosov and let $F$ be the vector field generating $\phi$.
Let $h\in C^{\infty}(M)$ and let $\theta$ be a smooth 1-form on $M$.
Then the cohomological equation
\[F(u)=h\circ \pi+\theta\]
has a solution $u\in C^{\infty}(SM)$ if and only if $h=0$ and $\theta$
is exact.}

\medskip

Note that by the smooth Liv\v sic theorem \cite{LMM} saying that
$h\circ\pi+\theta=F(u)$ is equivalent to saying that
$h\circ\pi+\theta$ has zero integral over every closed orbit of $\phi$.

Theorem B was proved in \cite{DP1} for the case of magnetic flows (i.e. $\la$ depends only on
$x$). It was surprising for us that the theorem also holds for systems that do not preserve
a smooth measure. The proof is also based on establishing a Pestov identity as 
in \cite{CS,DS} for geodesic flows, but some unexpected cancellations take place
producing in the end formulas which are just what one needs to prove the theorem.
%and we do not have as yet a conceptual explanation as to why exactly the proof works.
Earlier proofs of Theorem B for some geodesic and magnetic flows using Fourier analysis
can be found in \cite{GK,P2}.

Finally we note that Theorem A also holds if we allow magnetic forces. Indeed 
Theorem B holds for a generalized thermostat and $\Div F=-\theta$ even when
we have a magnetic field present. The extension of Theorem A to {\it isoenergetic}
thermostats (i.e. in the presence of potential forces) is discussed in Remark \ref{last}.

\bigskip

{\it Acknowledgements:} The first author thanks the Department of Pure Mathematics and
Mathematical Statistics at the University of Cambridge and Trinity College for hospitality 
and financial support while this work was in progress.

\section{Preliminaries}

Let $M$ be a closed oriented surface, $SM$ the unit sphere bundle
and $\pi:SM\to M$ the canonical projection. The latter is in fact
a principal $S^{1}$-fibration and we let $V$ be the infinitesimal
generator of the action of $S^1$.

Given a unit vector $v\in T_{x}M$, we will denote by $iv$ the
unique unit vector orthogonal to $v$ such that $\{v,iv\}$ is an
oriented basis of $T_{x}M$. There are two basic 1-forms $\alpha$
and $\beta$ on $SM$ which are defined by the formulas:
\[\alpha_{(x,v)}(\xi):=\langle d_{(x,v)}\pi(\xi),v\rangle;\]
\[\beta_{(x,v)}(\xi):=\langle d_{(x,v)}\pi(\xi),iv\rangle.\]
The form $\alpha$ is the canonical contact form of $SM$ whose Reeb vector
field is the geodesic vector field $X$.
The volume form $\alpha\wedge d\alpha$ gives rise to the Liouville measure
$d\mu$ of $SM$.

A basic theorem in 2-dimensional Riemannian geometry asserts that
there exists a~unique 1-form $\psi$ on $SM$ (the connection form)
such that $\psi(V)=1$ and
\begin{align}
& d\alpha=\psi\wedge \beta\label{riem1}\\ & d\beta=-\psi\wedge
\alpha\label{riem2}\\ & d\psi=-(K\circ\pi)\,\alpha\wedge\beta
\label{riem3}
\end{align}
where $K$ is the Gaussian curvature of $M$. In fact, the form
$\psi$ is given by
\[\psi_{(x,v)}(\xi)=\left\langle \frac{DZ}{dt}(0),iv\right\rangle,\]
where $Z:(-\varepsilon,\varepsilon)\to SM$ is any curve with
$Z(0)=(x,v)$ and $\dot{Z}(0)=\xi$ and $\frac{DZ}{dt}$ is the
covariant derivative of $Z$ along the curve $\pi\circ Z$.

%It is easy to check that $\alpha\wedge\beta=\pi^{*}\Omega_{a}$,
%hence
%\begin{equation}
%d\psi=-\pi^{*}(K\,\Omega_{a}). \label{psi}
%\end{equation}

For later use it is convenient to introduce the vector field $H$
uniquely defined by the conditions $\beta(H)=1$ and
$\alpha(H)=\psi(H)=0$. The vector fields $X,H$ and $V$ are dual to
$\alpha,\beta$ and $\psi$ and as a consequence of (\ref{riem1}--\ref{riem3}) they satisfy the commutation relations
\begin{equation}\label{comm}
[V,X]=H,\quad [V,H]=-X,\quad [X,H]=KV.
\end{equation}
Equations (\ref{riem1}--\ref{riem3}) also imply that the vector fields
$X,H$ and $V$ preserve the volume form $\alpha\wedge d\alpha$ and hence
the Liouville measure.

\section{An integral identity}

Henceforth $(M,g)$ is a closed oriented surface
and $X$, $H$, and $V$ are the same vector fields on $SM$
as in the previous section.

Let $\lambda$ be the smooth function on $SM$ given by (\ref{eqgt}), and let
$$
F=X+\lambda V
$$
be the generating vector field of the generalized thermostat.

From (\ref{comm}) we obtain:
$$
[V,F]=H+V(\la)V,\quad [V,H]=-F+\lambda V,\quad 
[F,H]=-\lambda F+(K-H(\lambda)+\lambda^2)V.
$$
%Note that
%$$
%H\lambda(x,v)=\langle \nabla \lambda(x), iv\rangle.
%$$

\begin{Lemma}[The Pestov identity]\label{pestov}
For every smooth function $u:SM\to\mathbb R$ we have
\begin{align*}
2Hu\cdot VFu
&=(Fu)^2+(Hu)^2-(K-H(\lambda)+\lambda^2)(Vu)^2\\
&+F(Hu\cdot Vu)+V(\la)Hu\cdot Vu-H(Fu\cdot Vu)+V(Fu\cdot Hu).
\end{align*}
\end{Lemma}

%\begin{Remark}
%{\rm A similar identity for the vector fields 
%$X$, $H_\lambda:=H+\lambda V$ and $V$
%was obtained in \cite[Lemma 2.1]{SU}}.
%\end{Remark}

\begin{proof}
Using the commutation formulas, we deduce:
\begin{align*}
2Hu\cdot VF&u-V(Hu\cdot Fu)\\
&=Hu\cdot VFu-VHu\cdot Fu\\
&=Hu\cdot(FVu+[V,F]u)
-Fu\cdot(HVu+[V,H]u)\\
&=Hu\cdot(F Vu+Hu+V(\la)Vu)
-Fu\cdot(HVu-Fu+\lambda Vu)\\
&=(Fu)^2+(Hu)^2
+(F Vu)(Hu)-(HVu)(Fu)
-\la Fu\cdot Vu+Hu\cdot V(\la)Vu\\
&=(Fu)^2+(Hu)^2
+F(Vu\cdot Hu)
-H(Vu\cdot Fu)
-[F,H]u\cdot Vu\\
&-\lambda Fu\cdot Vu+Hu\cdot V(\la)Vu\\
&=(Fu)^2+(Hu)^2
+F(Vu\cdot Hu)+V(\la)Hu\cdot Vu
-H(Vu\cdot Fu)\\
&-(K-H(\lambda)+\lambda^2)(Vu)^2
\end{align*}
which is equivalent to the Pestov identity.
\end{proof}

Now let $\Theta:=\alpha\wedge d\alpha$. This volume form generates the
Liouville measure $d\mu$.

\begin{Lemma} We have:
\begin{align}
L_{F}\Theta &=V(\la)\Theta;\label{lie1}\\
L_{H}\Theta &=0;\label{lie2}\\
L_{V}\Theta &=0.\label{lie3}
\end{align}

\label{lied}
\end{Lemma}

\begin{proof} Note that for any vector field $Y$, $L_{Y}\Theta=d(i_{Y}\Theta)$.
Since
$i_{V}\Theta=-\alpha\wedge\beta=-\pi^*\Omega_{a}$, where $\Omega_{a}$ is the area
form of $M$, we see that $L_{V}\Theta=0$.
Similarly, $L_{X}\Theta=L_{H}\Theta=0$.
Finally $L_{F}\Theta=L_{X}\Theta+L_{\la V}\Theta=d(i_{\la V}\Theta)=
V(\la)\Theta$.

\end{proof}

Below we will use the following consequence of Stokes theorem. Let $N$ be a closed
oriented manifold and $\Theta$ a volume form. Let $X$ be a vector
field on $N$ and $f:N\to\re$ a smooth function. Then
\begin{equation}
\int_{N}X(f)\,\Theta=-\int_{N}f\,L_{X}\Theta.
\label{bas}
\end{equation}

Integrating the Pestov identity over $SM$ against the Liouville measure $d\mu$,
and using (\ref{lie2}) and (\ref{lie3}) we obtain:

\begin{align*}
2\int_{SM} Hu\cdot VFu\,d\mu
&=\int_{SM}(Fu)^2\,d\mu+\int_{SM}(Hu)^2\,d\mu\\
&+\int_{SM}(F(Hu\cdot Vu)+V(\la)Hu\cdot Vu)\,d\mu\\
&-\int_{SM}(K-H(\lambda)+\lambda^2)(Vu)^2\,d\mu.\notag
\end{align*}
Using (\ref{bas}) and (\ref{lie1}) we get:
\[\int_{SM}(F(Hu\cdot Vu)+V(\la)Hu\cdot Vu)\,d\mu=0\]
and thus

\begin{align}\label{id1}
2\int_{SM} Hu\cdot VFu\,d\mu
&=\int_{SM}(Fu)^2\,d\mu+\int_{SM}(Hu)^2\,d\mu\\
&-\int_{SM}(K-H(\lambda)+\lambda^2)(Vu)^2\,d\mu.\notag
\end{align}

We will derive one more integral identity.
By the commutation relations, we have
$$
F Vu=VF u-Hu-V(\la)Vu.
$$
Therefore,
\[(F Vu)^2=(VF u)^2+(Hu)^2+(V(\la))^2(Vu)^2\]
\[-2VF u\cdot Hu-2VF u\cdot V(\la)Vu+2V(\la)Vu\cdot Hu.\]
Thus using again the commutation relations:
\[(F Vu)^2=(VF u)^2+(Hu)^2+(V(\la))^2(Vu)^2\]
\[-2VF u\cdot Hu-2F Vu\cdot V(\la)Vu-2(V(\la))^2(Vu)^2.\]
Since
\[F(V(\la)(Vu)^2)=2V(\la)Vu\cdot FVu+(Vu)^2F(V(\la))\]
we obtain:
\[(F Vu)^2=(VF u)^2+(Hu)^2-(V(\la))^2(Vu)^2\]
\[-2VF u\cdot Hu-F(V(\la)(Vu)^2)-(Vu)^2F(V(\la)).\]
Integrating this equation we obtain:
\begin{align}\label{id2}
2\int_{SM} Hu\cdot VFu\,d\mu
&=\int_{SM}(VFu)^2\,d\mu+\int_{SM}(Hu)^2\,d\mu\\
&-\int_{SM}(FVu)^2\,d\mu-\int_{SM}F(V(\la))(Vu)^2\,d\mu\notag
\end{align}
since by (\ref{bas}) and (\ref{lie1}) we get:
\[\int_{SM}\{F(V(\la)(Vu)^2)+(V(\la))^2(Vu)^2\}\,d\mu=0.\]

Combining (\ref{id1}) and (\ref{id2}) we arrive at the final integral identity
of this section:

\begin{Theorem}
\begin{equation}\label{id}
\int_{SM}(F Vu)^2\,d\mu-\int_{SM}{\mathbb K}(Vu)^2\,d\mu
=\int_{SM}(VF u)^2\,d\mu-\int_{SM}(Fu)^2\,d\mu,
\end{equation}
where ${\mathbb K}:=K-H(\la)+\la^2+F(V(\la))$.
\end{Theorem}

Of course this identity holds without any assumption on the underlying dynamics.
In the next section we will show how to use the Anosov hypothesis
to rewrite the left hand side of (\ref{id}) in terms of the stable or unstable
bundles.
At this point the proof differs from the one presented in \cite{DP1}.
We can no longer estimate the left hand side of (\ref{id}) using closed orbits
and the non-negative Liv\v sic theorem \cite{LT,PS} since in our context the Liouville
measure is not necessarily invariant.

\section{Using the Anosov property}

Recall that the Anosov property means that $T(SM)$
splits as $T(SM)=\re F\oplus E^{u}\oplus E^{s}$ in such a way that
there are constants $C>0$ and $0<\rho<1<\eta$ such that 
for all $t>0$ we have
\[\|d\phi_{-t}|_{E^{u}}\|\leq C\,\eta^{-t}\;\;\;\;\mbox{\rm
and}\;\;\;\|d\phi_{t}|_{E^{s}}\|\leq C\,\rho^{t}.\]
The subbundles are then invariant and H\"older continuous and have
smooth integral manifolds, the stable and unstable manifolds,
which define a continuous foliation with smooth leaves.

Let us introduce the weak stable and unstable bundles:
\[E^+=\re F\oplus E^s,\]
\[E^-=\re F\oplus E^u.\]

\begin{Lemma} For any $(x,v)\in SM$, $V(x,v)\notin E^{\pm}(x,v)$.
\end{Lemma}

\begin{proof} Let $\Lambda(SM)$ be the bundle over $SM$ such that at each point $(x,v)\in SM$
consists of all 2-dimensional subspaces $W$ of $T_{(x,v)}SM$ with
$F(x,v)\in W$.

The map $(x,v)\mapsto \V:=\re F(x,v)\oplus \re V(x,v)$ is a section of 
$\Lambda(SM)$ and its image is a codimension one submanifold that we denote by $\Lambda_{V}$.
Similarly the map $(x,v)\mapsto \re F(x,v)\oplus \re H(x,v)$ is a section of 
$\Lambda(SM)$ and its image is a codimension one submanifold that we denote by $\Lambda_{H}$.

The flow $\phi$ naturally lifts to a flow $\phi^*$ acting on $\Lambda(SM)$ via its differential.
Let $F^*$ be the infinitesimal generator of $\phi^*$.

{\bf Claim.} $F^*$ is transversal to $\Lambda_V$.

To prove the claim we define a function $m:\Lambda(SM)\setminus \Lambda_{H}\to\re$
as follows. If $W\in\Lambda(SM)\setminus \Lambda_{H}$, then $H\notin W$.
Thus there exists a unique $m=m(W)$ such that $mH+V\in W$. Clearly $m$ is smooth
and $\Lambda_{V}=m^{-1}(0)\subset\Lambda(SM)\setminus \Lambda_{H}$.
Fix $(x,v)\in SM$ and
set $m(t):=m(\phi^*_{t}(\V(x,v)))$. By the definition of $m$, there exist functions 
$x(t)$ and $y(t)$ such that

\[m(t)H(t)+V(t)=x(t)F(t)+y(t)d\phi_{t}(V).\]
Equivalently
\[m(t)d\phi_{-t}(H(t))+d\phi_{-t}(V(t))=x(t)F+y(t)V.\]
Differentiating with respect to $t$ and setting $t=0$ (recall that $m(0)=0$)
we obtain:
\[\dot{m}(0)H+[F,V]=\dot{x}(0)F+\dot{y}(0)V.\]
But $[V,F]=H+V(\la)V$. Thus $\dot{m}(0)=1$ which proves the Claim.

From the Claim it follows that $\Lambda_V$ determines an oriented
codimension one cycle in $\Lambda(SM)$ and by duality it defines
a cohomology class ${\mathfrak m}\in H^{1}(\Lambda(SM),\Z)$. Set $E=E^\pm$.
Given a continuous closed curve $\alpha:S^1\to SM$, the {\it index}
of $\alpha$ is $\nu(\alpha):=\langle {\mathfrak m}, [E\circ\alpha]\rangle$
(i.e. $\nu=E^*{\mathfrak m}\in H^{1}(SM,\Z)$).
The index of $\alpha$ only depends on the homology class
of $\alpha$.
Since $E$ is $\phi$-invariant, the Claim also ensures that if
$\ga$ is any closed orbit of $\phi$, then $\nu(\ga)\geq 0$.

Recall that according to Ghys \cite{Ghy} we know that $\phi$ is topologically
conjugate to the geodesic flow of a metric of constant negative curvature.
In particular, every homology class in $H_{1}(SM,Z)$ contains
a closed orbit of $\phi$. Thus $\nu$ must vanish.

If there exists $(x,v)\in SM$ for which
$V(x,v)\in E(x,v)$, then using that every point of $\phi$ is non-wandering,
we can produce exactly as in \cite[Lemma 2.49]{P1} a closed curve
$\alpha:S^1\to SM$ with $\nu(\alpha)>0$. 
This contradiction shows the lemma.

\end{proof}

\begin{Remark}{\rm The reader will recognize that the index that appears
in the proof of the lemma reduces to the Maslov index when $\phi$ is Hamiltonian.
The proof of the lemma also follows the presentation in \cite[Chapter 2]{P1} of analogous
results for geodesic flows.}
\end{Remark}

The lemma implies that there exist unique continuous functions $r^{\pm}$ on $SM$
such that
\[H+r^+ V\in E^+,\]
\[H+r^- V\in E^-.\]
Note that the Anosov property implies that $r^+\neq r^-$ everywhere.
Below we will need to use that the functions $r^{\pm}$ satisfy a Riccati
type equation along the flow. Note that $r^{\pm}$ are smooth along $\phi$
because $E^\pm$ are $\phi$-invariant.

\begin{Lemma} Let $r=r^{\pm}$. Then
\[F(r-V(\la))+r(r-V(\la))+\K=0.\]
\label{riccati}
\end{Lemma}

\begin{proof}Let $E=E^\pm$. Fix $(x,v)\in SM$, flow along $\phi$ and set
\[\xi(t):=d\phi_{-t}(H(t)+r(t)V(t)).\]
By the definition of $r$, $\xi(t)\in E(x,v)$ for all $t$.
Differentiating with respect to $t$ and setting $t=0$
we obtain:
\[\dot{\xi}(0)=[F,H]+F(r)V+r[F,V].\]
Using that
$$
[V,F]=H+V(\la)V,\quad 
[F,H]=-\lambda F+(K-H(\lambda)+\lambda^2)V
$$
we have
\[\dot{\xi}(0)=-\la F-rH+\left\{F(r)+K-H(\lambda)+\lambda^2-V(\la)r\right\}V.\]
Replacing $H$ by $\xi(0)-rV$ yields:
\[\dot{\xi}(0)+r\xi(0)+\la F=\left\{r^2+F(r)+K-H(\lambda)+\lambda^2-V(\la)r\right\}V.\]
Since $\dot{\xi}(0)+r\xi(0)+\la F\in E$ we must have
\[r^2+F(r)+K-H(\lambda)+\lambda^2-V(\la)r=0\]
which is the desired equation since
${\mathbb K}=K-H(\la)+\la^2+F(V(\la))$.

\end{proof}

Here is the main result of this section:

\begin{Theorem} Let $\psi:SM\to\re$ be a smooth function and suppose $\phi$
is Anosov. Then for $r=r^{\pm}$
\[\int_{SM}(F \psi)^2\,d\mu-\int_{SM}\K \psi^2\,d\mu=
\int_{SM}[F(\psi)-r\psi+\psi V(\la)]^2\,d\mu\geq 0.\]
Moreover,
\[\int_{SM}[F(\psi)-r\psi+\psi V(\la)]^2\,d\mu=0\]
if and only if $\psi=0$.
\label{llave2}
\end{Theorem}

\begin{proof} Let us expand $[F(\psi)-r\psi+\psi V(\la)]^2$:
\begin{align*}
[F(\psi)-r\psi+\psi V(\la)]^2&=[F(\psi)]^2+\psi^2 r^2+\psi^2[V(\la)]^2\\
&-2F(\psi)\psi r+2F(\psi)\psi V(\la)-2\psi^2 r V(\la).
\end{align*}
Using that (see Lemma \ref{riccati})
\[F(r-V(\la))+r(r-V(\la))+\K=0\]
we obtain:
\begin{align*}
[F(\psi)-r\psi+\psi V(\la)]^2&=[F(\psi)]^2-\K\psi^2\\
&-F((r-V(\la))\psi^2)+\psi^2[V(\la)]^2-\psi^2 r V(\la).
\end{align*}
If we integrate the last equality with respect to the Liouville measure $\mu$
we obtain as desired:
\[\int_{SM}(F \psi)^2\,d\mu-\int_{SM}\K \psi^2\,d\mu=
\int_{SM}[F(\psi)-r\psi+\psi V(\la)]^2\,d\mu\]
since by (\ref{bas}) and (\ref{lie1}) we have the following cancellation:
\[\int_{SM}\{-F((r-V(\la))\psi^2)+\psi^2[V(\la)]^2-\psi^2 r V(\la)\}\,d\mu=0.\]

Suppose now
\[\int_{SM}[F(\psi)-r\psi+\psi V(\la)]^2\,d\mu=0\]
which implies
\[F(\psi)-r\psi+\psi V(\la)=0\]
everywhere. Since this holds for $r=r^{\pm}$ we deduce:
\[(r^{+}-r^{-})\psi=0.\]
But for an Anosov flow $r^+-r^-\neq 0$ everywhere, thus $\psi=0$.

\end{proof}

\section{Proof of Theorem B}

Let us now prove Theorem B.
If $Fu=h\circ\pi+\theta$, then it is easy to see that
the right-hand side of (\ref{id}) is nonpositive.
Indeed, since $\mu$ is invariant under $v\mapsto -v$ and $v\to iv$ we have
\[\int_{SM}\theta_{x}(v)\,d\mu=0\;\;\; \mbox{\rm and}\;\; \int_{SM}(\theta_{x}(v))^2\,d\mu=
\int_{SM}(\theta_{x}(iv))^2\,d\mu.\]
But $VFu=\theta_{x}(iv)$ and thus
\[\int_{SM}(VFu)^2\,d\mu-\int_{SM}(Fu)^2\,d\mu=
-\int_{SM}(h\circ\pi)^{2}\,d\mu\leq 0.\]
Setting $\psi=Vu$, we get
\begin{equation}\label{in}
\int_{SM}\left\{(F \psi)^2
-\K \psi^2\right\}\,d\mu\le 0.
\end{equation}
By Theorem \ref{llave2} this happens if and only if $\psi=0$.
This would give $Vu=0$, which says that $u=f\circ\pi$
where $f$ is a smooth function on $M$. But in this case, since
$d\pi_{(x,v)}(F)=v$ we have
$Fu=df_{x}(v)$. This clearly implies the claim of the theorem.

\begin{Remark}{\rm Suppose that we include potential forces in our dynamics, that is,
we consider the {\it isoenergetic thermostat}:

\begin{equation}
\frac{D\dot{\gamma}}{dt}=-\nabla W+E(\ga)-\frac{\langle E(\ga),\dot{\ga}\rangle}{|\dot{\ga}|^2}\,\dot{\ga}
\label{eqt1}
\end{equation}
on the energy level $\frac{1}{2}|v|^2+W(x)=k$ (we assume that $|v|$ does not vanish on the energy level).
Wojtkowski has pointed out \cite[Theorem 2.4]{W3} that the dynamics of (\ref{eqt1}) reparametrized
by arc-length defines a flow on $SM$ which coincides with the isokinetic thermostat with external field
\[\widetilde{E}:=\frac{-\nabla W+E}{2(k-W)}=\frac{1}{2}\nabla(\log(k-W))+\frac{E}{2(k-W)}.\]
Since the vanishing of entropy production and the Anosov property are unaltered by smooth time changes
we conclude applying Theorem A to $\widetilde{E}$ that an Anosov isoenergetic thermostat has zero entropy
production if and only if $E/2(k-W)$ has a global potential.

The question of whether Theorem B extends to higher dimension is more delicate. We hope to discuss
this topic elsewhere.

}
\label{last}
\end{Remark}

\end{document}